\documentclass{amsart}
\def\data{05 February 2024}
\usepackage{amssymb}
\usepackage{amsfonts}
\usepackage{amsmath}
\usepackage{amsthm}
\usepackage[pdftex,bookmarks,colorlinks]{hyperref}
\usepackage[pdftex]{hyperref}
\usepackage{url}
\usepackage{color}
\usepackage{graphicx}

\newcommand{\bq}{\begin{quote}}
\newcommand{\eq}{\end{quote}}
\newcommand{\bi}{\begin{itemize}}
\newcommand{\ei}{\end{itemize}}
\newcommand{\bd}{\begin{description}}
\newcommand{\ed}{\end{description}}
\newcommand{\ben}{\begin{enumerate}}
\newcommand{\een}{\end{enumerate}}
\newcommand{\bbm}{\begin{bmatrix}}
\newcommand{\ebm}{\end{bmatrix}}
\newcommand{\bea}{\begin{eqnarray*}}
\newcommand{\eea}{\end{eqnarray*}}

\newtheorem{theorem}{Theorem}

\def\bF{\mathbb{F}}

\def\lsl{\mathfrak{sl}}
\def\lgl{\mathfrak{gl}}

\def\liz{\mathfrak{z}}

\def\CC{\mathsf{C}}
\def\DD{\mathsf{D}}
\def\EE{\mathsf{E}}

\def\KK{\mathsf{K}}

\def\NN{\mathsf{N}}

\def\RR{\mathsf{R}}
\def\SS{\mathsf{S}}
\def\TT{\mathsf{T}}

\def\XX{\mathsf{X}}
\def\YY{\mathsf{Y}}
\def\ZZ{\mathsf{Z}}

\newcommand{\bx}[1]{\mathsf{#1}}     

\def\2G2{\ensuremath{^2{\rm G}_2}}
\def\sl{{\rm{SL}}}
\def\gl{{\rm{GL}}}

\def\psl{{\rm{PSL}}}
\def\pgl{{\rm{PGL}}}

\def\so{{\rm{SO}}}

\definecolor{darkgreen}{rgb}{0,0.6,0}
\newcommand{\avb}[1]{{\color{darkgreen}#1}}

\newcommand{\encr}{\ensuremath{\vDash}}

\AtBeginDocument{%
  \def\MR#1{}
  }

\makeatletter
\@namedef{subjclassname@2010}{$2020$ Mathematics Subject Classification}
\makeatother

\begin{document}

\title[Black box rings of 2 by 2 matrices]{Structural proxies for black box rings encrypting rings of 2 by 2 matrices over finite fields of odd order}
\date{\data}

\author{Alexandre Borovik}
\address{Department of Mathematics, University of Manchester, UK\\ alexandre$\gg{\rm at}\ll$borovik.net}
\author{\c{S}\"{u}kr\"{u} Yal\c{c}\i nkaya}
\address{Department of Mathematics, Istinye University, Turkey\\  sukru.yalcinkaya$\gg{\rm at}\ll$istinye.edu.tr}


\begin{abstract} This paper provides an example of \emph{structural proxies} in the sense of our papers \cite{BY2017C,BY2020-sl2} for a class of algebraic structure mentioned in \cite{BY2017C}: black box rings encrypting rings of $2\times 2$ matrices of finite fields of odd order.
\end{abstract}

\maketitle

\section{Introduction}

General discussion of problems in black box algebra, and extensive bibliorgraphy and historic remarks on black box groups, fields and rings can be found in \cite{BY2018}; terminology and notation follows more updated papers  \cite{BY2017C,BY2020-sl2}.

In this paper, we prove the following

\begin{theorem}\label{th:m22}
Let $\RR \encr {\rm M}_{2\times 2}(\bF)$ be a black box ring encrypting the ring of $2\times 2$ matrices over a finite field\/ $\bF$ of odd order. Assume that we know this order, $|\bF| = q$. Then we can construct, in time polynomial in $\log q$, a black box field $\ZZ \encr \bF$, and polynomial {\rm (}in $\log q${\rm )} time isomorphisms
\[
\RR \leftrightarrows {\rm M}_{2\times 2}(\ZZ).
\]
\end{theorem}

\section{Black box algebra}

\subsection{Axiomatic description of black box algebraic structures}
\label{sec:BB1--BB3}

A \emph{black box algebraic structure} $\XX$ is a black box (or an oracle, or a device, or an algorithm) operating with $0$--$1$ strings of uniform length which encrypt (not necessarily in a unique way) elements of some algebraic structure $A$: if $\bx{x}$ is one of these strings  then it corresponds to a unique (but unknown to us) element $\pi(\bx{x}) \in A$. We call the elements of $\XX$ \emph{cryptoelements}.

Our axioms for black boxes are the same as in \cite{BY2014,BY2018,BY2020-sl2}, but stated in a more formal language.
\bi
\item[\textbf{BB1}] On request, $\XX$ produces a `random' cryptoelement $\bx{x}$ as a string of fixed length $l(\XX)$, which depends on $\XX$, which encrypts an element $\pi(x)$ of some fixed explicitly given  algebraic structure $A$; this is done in  time polynomial in $l(\XX)$. When this procedure is repeated, the  elements $\pi(\bx{x}_1), \pi(\bx{x}_2),\dots$ are independent and uniformly distributed in $A$.
\ei

To avoid messy notation, we assume that operations on $A$ are unary or binary; a general case can be treated in exactly the same way.

\bi
\item[\textbf{BB2}] On request,  $\XX$ performs algebraic operations on the encrypted strings which correspond to operations in  $A$ in a way which makes the map $\pi$ (unknown to us!) a homomorphism: for every binary (unary case is similar) operation $\boxdot $ and strings $\bx{x}$ and $\bx{y}$ produced or computed by $\XX$,
\[
\pi(\bx{x} \boxdot \bx{y}) = \pi(\bx{x}) \boxdot \pi(\bx{y}).
\]
\ei

It should be noted that we do not assume the existence of an algorithm which allows us to decide whether a specific string can be potentially produced by $\XX$; requests for operations on strings can be made only in relation to cryptoelements previously output by $\XX$. Also, we do not make any assumptions on probabilistic distribution of cryptoelements.

\begin{itemize}
\item[\textbf{BB3}] On request, $\XX$ determines, in time polynomial in $l(\XX)$, whether two cryptoelements $\bx{x}$ and $\bx{y}$ encrypt the same element in $A$, that is, check whether $\pi(\bx{x}) = \pi(\bx{y})$.
\end{itemize}

We say in this situation that  a black box $\XX$ \emph{encrypts} the algebraic structure $A$; we shall denote that $\XX \encr A$. We use the same notation for elements: if $\pi(\bx{x}) = a$ for $\bx{x} \in \XX$ and $a \in A$, we also write $\bx{x} \encr a$, it could be convenient in calculations.

In our algorithms, we have to build new black boxes from existing ones and work with several black box structures at once: this is why we have to keep track of the length $l(\bx{\XX})$ on which a specific black box $\XX$ operates. For example, it turns out in \cite{BY2018} that it is useful to consider an automorphism of $A$ as a graph in $A\times A$. This produces an another algebraic structure isomorphic to $A$ which can be seen as being encrypted by a black box $\ZZ$ producing, and operating on, certain pairs of strings from $\XX$, see \cite{BY2018} for more examples. In this case, clearly, $l(\ZZ)  =2l(\XX)$.

We note that when we build a new black box from an existing black box algebraic structure, we do not produce a list of elements in this new black box algebraic structure but only algorithms which perform the tasks in the Axioms BB1, BB2 and BB3.

\subsection{Structural proxy}\label{sec:structural_proxy}

Most groups of Lie type (we exclude $^2B_2$, $^3D_4$, $^2F_4$ and $^2G_2$ to avoid technical details) can be seen as functors $G: \mathcal{F} \longrightarrow \mathcal{G}$ from the category of fields $\mathcal{F}$ with an automorphism of order $\leqslant 2$ to the category of groups $\mathcal{G}$. There are also other algebraic structures which can be defined in a similar way as functors from $\mathcal{F}$, for example projective planes or simple Lie algebras (viewed as rings). The following problem is natural and, as our results show, useful in this context.

\bi
\item[] \emph{Construction of a structural proxy:}  Suppose that we are given  a black box structure $\XX\encr A(\bF)$. Construct, in time polynomial  in $l(\XX)$,
    \bi
     \item a black box field $\KK\encr \bF$, and
     \item  two way bijective morphisms $A(\KK) \longleftrightarrow \XX$.
     \ei
\ei

If we construct a black box field $\KK$ by using $\XX$ as a computational engine, then we can construct the natural representation $A(\KK)$ of the structure $A$ over the black box field $\KK$. By Maurer and Raub  \cite{maurer07.427} (see also \cite[Theorem 3.1]{BY2018}), we have a computable in probabilistic polynomial time isomorphism $\bF_q \longrightarrow \KK$  and hence an isomorphism $A(\bF_q) \longrightarrow A(\KK)$ completing the structure recovery of $\XX$.

Structural proxies  and structure recovery play a crucial role in our papers \cite{BY2018,BY2020-sl2}. We construct there structural proxies for
\bi
\item a black box projective plane with polarity $\ZZ \encr \mathbb{P}(\bF_q)$;
\item the projectivisation of a black box Lie algebra $\YY\encr  \lsl_2(\bF_q)$ (this one appears in disguise as construction of the ``cross product'' on the projective plane);
\item black box groups $\XX_1\encr \so_3(\bF_q)\simeq \pgl_2(\bF_q)$, $\XX_2 \encr \psl_2(\bF_q)$, $\XX_3 \encr \sl_2(\bF_q)$.
\ei
In all these cases, $q$ is odd.

Some striking parallels with model-theoretic algebra are discussed in our paper \cite{BY2020-models}.

\section{Proof of Theorem \ref{th:m22}}

We start with an obvious observation that even if the zero $\bx{0} \encr 0$ is not given explicitly in the set-up of the black box ring $\RR$, we can  easily construct it by taking a random $\bx{r} \in \RR$ and computing  $\bx{0} = \bx{r}-\bx{r}$.

Construction of a structural proxy of $\RR$ is explained in the steps below. The crucial ingredient of the proof is to construct the field of scalars and a dihedral group of order 8 in the multiplicative group of $\RR$.

\medskip\noindent\textbf{Step 1: Construction of the black box group $\XX \encr \gl_2(\bF_q)$.}
Observe first that a random element in $R= {\rm M}_{2\times 2}(\bF_q)$ is invertible with probability $1 - O\left(\frac{1}{q}\right)$.  Therefore a random cryptoelement $\bx{r}\in \RR$ is invertible with probability close to $1$ when $q$ is large. Set
\[
E = q(q^2-1);
\]
then, for an invertible cryptoelement $\bx{r} \in \RR$, we have $\bx{r}^{E} = \bx{e}$ where $\bx{e}$ encrypts the multiplicative identity of $\RR$. Therefore we can easily construct the identity element in $\RR$. Hence  the invertible cryptoelements in $\RR$ give us the black box group $\XX \encr \gl_2(\bF_q)$.

\medskip\noindent\textbf{Step 2: Generation of trace $0$ elements.}
Our next task is to construct in $\RR$ a black box subring $\ZZ$ encrypting the subring $Z = Z(R)$ of scalar matrices in $R$.  Of course, $Z$ is a field, and therefore $\ZZ \encr \bF_q$. When $\ZZ$ is constructed, it will turn $\RR$ into a black box $\ZZ$-algebra.

It is easy to check that, for any matrices $a,b \in {\rm M}_{2\times 2}(\bF_q)$, the matrix $[a,b]^2 = (ab-ba)^2$ is a scalar matrix; hence we have a natural map
\bea
\RR \times \RR & \longrightarrow & \ZZ\\
(\bx{r},\bx{s}) & \mapsto & [\bx{r},\bx{s}]^2
\eea
with values in $\ZZ = Z(\RR)$. We need to check that this map gives almost uniformly distributed cryptoelements in $\ZZ$.

Matrices of trace zero in ${\rm M}_{2\times 2}(\bF_q)$ form the Lie subalgebra $\lsl_2(\bF_q)$. The Lie algebra $\lgl_2(\bF_q)$ of ${\rm M}_{2\times 2}(\bF_q)$ can be decomposed as the sum of two ideals
\[
\liz \oplus \lsl_2(\bF_q),
\]
where $\liz$ is the center of $\lgl_2(\bF_q)$. Therefore the probability of an element $c \in \lsl_2(\bF_q)$ to be a commutator of two independent random elements from $\lgl_2(\bF_q)$ is the same as the probability of being a commutator of two independent random elements from $\lsl_2(\bF_q)$. For estimating the latter, observe that if $c = [a,b]$ and $c \ne 0$, then $a$ and $b$ belong to the plane in $\lgl_2(\bF_q)$ orthogonal to $c$ {with respect to the Killing form on $\lgl_2(\bF_q)$; this plane does not contain $c$, if $c$ is a semisimple element, and contains $c$, if $c$ is nilpotent. It is easy to see that probability of $c$ being a commutator of random elements from $\lsl_2(\bF_q)$ is $\frac{1}{q^3} +O\left(\frac{1}{q^4}\right)$.

\medskip\noindent\textbf{Step 3: The quadratic form on $\lsl_2$.}
Now we turn our attention to the square map
\bea
\sigma: \lsl_2(\bF_q) & \longrightarrow & Z\\
 s & \mapsto & s^2.
\eea
Let $s$ be a trace zero matrix,
\[
s = \bbm a & b \\ c & -a \ebm,
\]
then it is easy to check that
\[
s^2 = \bbm  a^2 +bc  & 0 \\ 0 &  a^2 +bc \ebm,
\]
and that $a^2 + bc = -\det s$
is proportional to the Killing quadratic form on $\lsl_2(\bF_q)$. It is easy to see that elements of $Z$ are images of random elements from $\lsl_2(\bF_q)$ with probability $\frac{1}{q} + O\left(\frac{1}{q^2}\right)$.

\medskip\noindent\textbf{Step 4: Construction of the field of scalars.}
Combining these two esimates, we see that the square map
\bea
R \times R & \longrightarrow & Z\\
(a,b) & \mapsto & [a,b]^2\\
\eea
hits specific elements in $Z$ with probabilities $\frac{1}{q} + O\left(\frac{1}{q^2}\right)$, that is, it is essentially the uniform distribution. Therefore the map
\bea
\RR \times \RR & \longrightarrow & \ZZ\\
(\bx{x},\bx{y}) & \mapsto & [\bx{x},\bx{y}]^2\\
\eea
can be taken for a generator of random elements for the black box $\ZZ$; operations of addition, multiplication, inversion on $\ZZ$ are inherited from $\RR$ and its multiplicative group.

\medskip\noindent\textbf{Step 5: The normaliser of a maximal torus in $\XX$.} We start by constructing a non-central involution $\bx{t} \in \XX$; it is done in the same way as it has been done in the black box group encrypting $\pgl_2(\bF_q)$ in \cite{BY2018}. Since all non-central involutions in the group $\gl_2(\bF_q)$ are conjugate, we can assume without loss of generality that $\bx{e}_1$ encrypts the matrix $\bbm -1 & 0 \\ 0 & 1 \ebm$, and we write
\[
\bx{e}_1 \encr \bbm -1 & 0 \\ 0 & 1 \ebm.
\]
The centraliser $\CC_{\XX}(\bx{e}_1)$ is the maximal torus
\[
\TT \encr T =\left\{\bbm * & 0 \\ 0 & * \ebm\right\},
\]
it is constructed by the standard algorithm for the construction of centralisers of involutions in black box groups. Then, by constructing involutions from random elements in $\TT$, we construct an} elementary abelian group $\EE$ of order $4$ in $\TT$. There is only one conjugacy class of elementary abelian subgroups of order $4$ in $\gl_2(\bF_q)$. Therefore, without loss of generality, we may assume that the elements of $\EE$ are
\[
\bx{e} \encr \bbm 1 & 0 \\ 0 & 1\ebm, \quad
-\bx{e} \encr \bbm -1 & 0 \\ 0 & -1\ebm, \quad
\bx{e}_1 \encr \bbm -1 & 0 \\ 0 & 1\ebm, \quad
\bx{e}_2 \encr \bbm 1 & 0 \\ 0 & -1\ebm.
\]
Next we take the factor group $\tilde{\XX} = \XX / \ZZ \encr \pgl_2(\bF_q)$ by replacing the equality relation in $\XX$ by the new one:
\[
\bx{x} \equiv \bx{y}
\]
if and only if $\bx{x}^{-1}\bx{y}$ commutes with elements $\bx{r},\bx{s}$  of odd order  and $[\bx{r},\bx{s}] \neq \bx{1}$. Then we construct, in a similar way, the image $\tilde{\TT}$ of $\TT$ in $\tilde{\XX}$; the corresponding subgroup $\tilde{T}$ in $X$ is a torus of order $q-1$ and therefore contains a unique involution $\tilde{t}$; we construct $\tilde{\bx{t}} \in \tilde{\TT}$ and its centraliser $\tilde{\NN}=\CC_{\tilde{\XX}}(\tilde{\bx{t}}) = \NN_{\tilde{\XX}}(\tilde{\TT})$. The torus $\tilde{\TT}$ has index $2$ in $\tilde{\NN}$; a random cryptoelement in $\tilde{\NN}$ belongs to $\tilde{\NN} \smallsetminus \tilde{\TT}$ with probability $\frac{1}{2}$; so we can pick $\tilde{\bx{w}} \in \tilde{\NN} \smallsetminus \tilde{\TT}$; it has order $2$ in $\tilde{\XX}$. Let $\bx{w}$ be its preimage in $\XX$, where we will be working from this point on. Now $\NN= \NN_{\XX}(\TT)=\TT\langle \bx{w} \rangle$ is the normaliser of a maximal torus in $\XX$.

\medskip\noindent\textbf{Step 6: A Sylow $2$-subgroup in  $\NN$.} Factorise $q-1 = 2^kl$ with $l$ odd. The map $\bx{t} \mapsto \bx{t}^l$ is a homomorphism from $\TT$ onto its Sylow $2$-subgroup  \avb{$\SS \encr Z_{2^k}\times Z_{2^k}$} (the \avb{direct} product of two cyclic subgroups of order $2^k$.

Compute $\bx{r} = \bx{w}^l$, then $\bx{r}$ is a $2$-element and $\SS\langle \bx{r} \rangle \encr Z_{2^k}\wr Z_2$ because $\NN$ contains an element encrypting the involution $\bbm 0 & 1 \\ 1 & 0\ebm$. We need to construct a dihedral group of order 8 in $\SS\langle \bx{r} \rangle$ and to do that we first construct a sequence of subgroups
\[
\bx{1} = \SS_0 <\EE = \SS_1 < \SS_2 < \dots < \SS_{k-1} < \SS_k = \SS
\]
such that $\SS_{i+1}/\SS_i \encr Z_2 \times Z_2$ for all $i=0,1,\dots,k-1$. We start by constructing $\SS$ by finding two generators $\bx{s}_1,\bx{s}_2$ of $\SS$. Then, we have
\[
\SS_{k-i} = \langle \bx{s}_1^{2^i}, \bx{s}_2^{2^i}\rangle, \, i=0,1, \ldots, k-1.
\]
To construct $\bx{s}_1,\bx{s}_2$, pick two random elements $\bx{t}_1,\bx{t_2}\in \TT$ and compute $\bx{s}_1 = \bx{t}_1^l$, $\bx{s}_2 = \bx{t}_2^l$. Observe that $\langle \bx{s}_1,\bx{s}_2 \rangle =\SS$ if and only if $\langle \bx{s}_1^{2^{k-1}},\bx{s}_2^{2^{k-1}}\rangle = \EE$. Now observe that $\bx{s}_1$ and $\bx{s}_2$ generate $\SS$ precisely when they both belong to $\SS\smallsetminus \SS_{k-1}$ which happens with probability
\[
\frac{4-1}{4} \times \frac{4-1}{4} = \frac{9}{16},
\]
and {they get into different cosets of $\SS$ over $\SS_{k-1}$ with probability \[
1 - 3 \cdot\frac{1}{3}\cdot\frac{1}{3} = \frac{2}{3}.
\]
Hence $\bx{s}_1$ and $\bx{s}_2$ generate $\SS$ with probability
\[
\frac{9}{16} \cdot \frac{2}{3} = \frac{3}{8}.
\]
After $n$ random tests, we find the generators with probability
\[
1- \left(\frac{5}{8}\right)^n.
\]

Observe that the criterion for an element $\bx{s}\in \SS = \SS_k$ to belong to $\SS_i$ is very simple: $\bx{s}^{2^i} = \bx{e}$. This allows us to compute in factor groups   $\SS\langle \bx{r}\rangle /\SS_i$, $i = 0,1,\dots, k-1$.

\medskip\noindent\textbf{Step 7: A dihedral subgroup of order $\boldsymbol{8}$ in $\XX$.}
Starting from $i=k-1$ and going down to $i=0$, we construct dihedral subgroups $\DD_{i}$ in $\SS\langle \bx{r}\rangle/\SS_i$ recursively. Denote $\bx{r}_k = \bx{r}$.  If $i=k-1$, we set $\DD_{k-1}= \SS_k\langle \bx{r}_k\rangle/\SS_{k-1}$, it is a non-abelian group of order $8$ which contains an elementary abelian group $\SS_k/\SS_{k-1}$, that is, $\DD_{k-1}$ is the dihedral group of order $8$. We identify, by direct inspection, an involution in $\DD_{k-1}$ and denote its preimage in $\SS_k\langle \bx{r}_k \rangle$ by $\bx{r}_{k-1}$. Then we consider the factor group $\DD_{k-2} = \SS_{k-1}\langle \bx{r}_{k-1} \rangle/ \SS_{k-2}$ and repeat the process until we get to the desired dihedral group $\DD_0$.

\medskip\noindent\textbf{Step 8: Matrix units in $\RR$.} As it can be easily seen from the character table of the  dihedral group of order $8$, all such subgroups in $\gl_2(\bF_q)$ are conjugate in $\gl_2(\bF_q)$. We retain notation for elements in $\EE < \DD_0$; we can assume, without loss of generality, that the involution $\bx{r}_1\in \DD_0$ found at the previous step represents the matrix $\bbm 0 & 1 \\ 1& 0\ebm$.

We compute
\[
\bx{e}_{11} = \frac{\bx{1}}{\bx{2}}\cdot (\bx{e} +\bx{e}_2) \; \mbox{ and } \;  \bx{e}_{22} = \frac{\bx{1}}{\bx{2}}\cdot (\bx{e} +\bx{e}_1),
\]
then
\[
\bx{e}_{11}\encr  \bbm 1 & 0 \\ 0 & 0 \ebm = e_{11} \; \mbox{ and } \;
\bx{e}_{22}\encr  \bbm 0 & 0 \\ 0 & 1 \ebm = e_{22}.
\]
Moreover,
\[
\bx{e}_{21} = \bx{r}_{1}\bx{e}_{11} \encr \bbm 0 & 1 \\ 1 & 0\ebm \bbm 1 & 0 \\ 0 & 0 \ebm  =\bbm 0 & 0 \\ 1 & 0 \ebm =  e_{21}
\]
and
\[
\bx{e}_{12} = \bx{r}_1\bx{e}_{22} \encr \bbm 0 & 1 \\ 1 & 0\ebm \bbm 0 & 0 \\ 0 & 1 \ebm = \bbm 0 & 1 \\ 0 & 0 \ebm = e_{12}.
\]

\medskip\noindent\textbf{Step 9: Representation of cryptoelements in $\RR$ by $2\times 2$ matrices over $\ZZ$.}  Take an arbitrary   $\bx{x} \in \RR$ and assume that
\[
\bx{x} \encr\bbm a_{11} & a_{12}\\ a_{21} & a_{22} \ebm, \,  a_{ij} \in \bF_q.
\]
We should construct proxies for $a_{ij}$, $i,j=1,2$, that is, we should construct elements $\bx{z}_{ij}\in \ZZ$ such that $\bx{z}_{ij} \encr a_{ij}$ and}\[
\bx{x} = \sum_{i,j} \bx{z}_{ij}\bx{e}_{ij}.
\]
To do that, we start by computing cryptoelements $\bx{x}_{ij} \in \RR$:
\bea
\bx{x}_{11} &:=& \bx{e}_{11}\bx{x}\bx{e}_{11} \encr a_{11}e_{11}\\
\bx{x}_{22} &:=& \bx{e}_{22}\bx{x}\bx{e}_{22} \encr a_{22}e_{22}\\
\bx{x}_{12} &:=& \bx{e}_{11}\bx{x}\bx{e}_{22} \encr a_{12}e_{12}\\
\bx{x}_{21} &:=& \bx{e}_{22}\bx{x}\bx{e}_{11} \encr a_{21}e_{21}\\
\eea
Then we produce scalar matrices:
\[
\bx{z}_{11} := \bx{x}_{11} + \bx{x}_{11}^{\bx{r}_1} \encr \bbm a_{11} & 0\\ 0 & a_{11} \ebm,
\]
hence
\[
\bx{x}_{11} = \bx{z}_{11}\bx{e}_{11}.
\]
Similarly,
\[
\bx{z}_{22} := \bx{x}_{22} + \bx{x}_{22}^{\bx{r}_1} \encr \bbm a_{22} & 0\\ 0 & a_{22} \ebm,
\]
hence
\[
\bx{x}_{22} = \bx{z}_{22}\bx{e}_{22}.
\]
Two other matrix element of $\bx{x}$ can be obtained as follows:
\[
\bx{z}_{12}:= \bx{r}_1\bx{x}_{12} + \bx{x}_{12}\bx{r}_1 \encr \bbm a_{12} & 0 \\ 0 & a_{12} \ebm,
\]
yielding
\[
\bx{x}_{12} = \bx{z}_{12}\bx{e}_{12},
\]
and similarly
\[
\bx{z}_{21}:= \bx{r}_1\bx{x}_{21} + \bx{x}_{21}\bx{r}_1 \encr \bbm a_{21} & 0 \\ 0 & a_{21} \ebm,
\]
yielding
\[
\bx{x}_{21} = \bx{z}_{21}\bx{e}_{21}.
\]
Hence, we have
\[
\sum_{i,j} \bx{z}_{ij}\bx{e}_{ij} \encr \bbm a_{11} & a_{12}\\ a_{21} & a_{22} \ebm,
\]
and
\[
\sum_{i,j} \bx{z}_{ij}\bx{e}_{ij} = \bx{x}
\]
which establishes a two way isomorphism between $\RR$ and $M_{2\times 2}(\ZZ)$.

\section*{Acknowledgements}

This paper---and other papers in our project---would have never been written if the authors did not enjoy the warm hospitality offered to them at the Nesin Mathematics Village in \c{S}irince, Izmir Province, Turkey,  as part of their Research in Pairs programme; our thanks go to Ali Nesin and to all volunteers, staff, and students who have made the Village a mathematical paradise.

\bibliographystyle{amsplain}

\providecommand{\bysame}{\leavevmode\hbox to3em{\hrulefill}\thinspace}
\providecommand{\MR}{\relax\ifhmode\unskip\space\fi MR }
\providecommand{\MRhref}[2]{%
  \href{http://www.ams.org/mathscinet-getitem?mr=#1}{#2}
}
\providecommand{\href}[2]{#2}

\end{document}